\documentclass[11pt,oneside]{amsart}
\usepackage{amsmath,amsthm,amssymb,amscd,fancybox}
\usepackage[a4paper,left=3cm,right=3cm,top=3cm,bottom=3cm]{geometry}
\usepackage{cite}
\renewcommand{\epsilon}{\varepsilon}
\renewcommand{\phi}{\varphi}
 \newcommand{\Lp}{L_{\vec{p}}(\bR^n)}
\newcommand{\Lpx}{L_{\vec{p}}(dx)}
 
 \newcommand{\bZ}{\mathbb{Z}}
\newcommand{\bR}{\mathbb{R}} \newcommand{\bN}{\mathbb{N}}
\newcommand{\cF}{\mathcal{F}} 
\newcommand{\dS}{\mathcal{S}(\bR^n)} 

\newcommand{\supp}{\operatorname{supp}}
\newcommand{\qp}{\mathcal{Q}}
\newcommand{\Zd}{\mathbb{Z}^n\setminus \{0\}}

\newcommand{\brac}[1]{\langle #1\rangle}
\newtheorem{theorem}{Theorem}[section]
\newtheorem{lemma}[theorem]{Lemma}
\newtheorem{proposition}[theorem]{Proposition}

\theoremstyle{definition}
\newtheorem{definition}[theorem]{Definition}
\newtheorem{example}[theorem]{Example}
\theoremstyle{remark}
\newtheorem{remark}[theorem]{Remark}
\numberwithin{equation}{section}

\def\bZ{{\mathbb Z}}

\def\bN{{\mathbb N}}

\def\bR{{\mathbb R}}

\def\cF{\mathcal{F}}

\def\cM{\mathcal{M}}

\begin{document}%\nocite{*}
\title[Pseudodifferential operators on Mixed-Norm $\alpha$-modulation spaces]{Pseudodifferential operators on Mixed-Norm $\alpha$-modulation spaces}
\author[M.\ Nielsen]{Morten Nielsen}
\address{Department of Mathematical Sciences\\ Aalborg
  University\\ Skjernvej 4A\\ DK-9220 Aalborg East\\ Denmark}
\email{mnielsen@math.aau.dk}
\begin{abstract}
  Mixed-norm $\alpha$-modulation spaces were introduced recently by Cleanthous and Georgiadis [Trans.\ Amer.\ Math.\ Soc.\ 373 (2020), no. 5, 3323-3356]. The mixed-norm spaces $M^{s,\alpha}_{\vec{p},q}(\bR^n)$,
  $\alpha\in [0,1]$, form a family of smoothness spaces that contain the 
  mixed-norm Besov spaces as special cases. In this paper we prove that
  a pseudodifferential operator $\sigma(x,D)$ with symbol in the
  H\"ormander class $S^b_{\rho}$ extends to a bounded operator
  $\sigma(x,D)\colon M^{s,\alpha}_{\vec{p},q}(\bR^n) \rightarrow
  M^{s-b,\alpha}_{\vec{p},q}(\bR^n)$ provided $0<\alpha\leq \rho\leq 1$,
  $\vec{p}\in (0,\infty)^n$, and $0<q<\infty$. The result extends the known result that
  pseudodifferential operators with symbol in the class $S^b_{1}$
  maps the mixed-norm Besov space $B^s_{\vec{p},q}(\bR^n)$ into $B^{s-b}_{\vec{p},q}(\bR^n)$.
\end{abstract}
\maketitle

\section{Introduction}

The $\alpha$-modulation spaces is a family of smoothness spaces that contains the Besov spaces and the
modulation spaces, introduced by Feichtinger \cite{fei}, as special
cases. In the non-mixed-norm setting, the $\alpha$-modulation spaces were introduced by Gr\"obner
\cite{Groebner1992}.   Gr\"obner used the general framework of decomposition type
Banach spaces considered by Feichtinger and Gr\"obner in
\cite{MR87b:46020,MR89a:46053} to build the $\alpha$-modulation
spaces.  The parameter $\alpha$ determines a specific type of
decomposition of the frequency space $\bR^n$ used to define
the space $M^{s,\alpha}_{p,q}(\bR^n)$.
The $\alpha$-modulation spaces contain the Besov spaces and the
modulation spaces, introduced by Feichtinger \cite{fei}, as special
cases. The choice $\alpha=0$ corresponds to the
classical modulation spaces $M^s_{p,q}(\bR^n)$, and $\alpha=1$ 
corresponds to the Besov scale of spaces.

 The applicability of $\alpha$-modulation spaces to the study of pseudo-differential operators comes rather natural, and in fact the family of coverings used to construct 
the $\alpha$-modulation spaces was considered independently by
P{\"a}iv{\"a}rinta and Somersalo in
\cite{MR975206}. P{\"a}iv{\"a}rinta and Somersalo used the partitions
to extend the  {C}alder\'on-{V}aillancourt boundedness result for
pseudodifferential operators to local Hardy spaces.

Recently, function spaces in anisotropic and mixed-norm settings have attached considerable interest, see for example \cite{AI, BB, MR2423282, Bownik2006,Bownik2005, MR3622656,MR3573690, JS} and reference therein. This is in part driven by advances in the study of partial and pseudodifferential operators, where there is a natural desire to be able to better model and analyse anisotropic phenomena.
In particular, pseudo-differential operators on mixed-norm Besov spaces have been studied in  \cite{MR3682609,Cleanthous2019}.

In this paper we study pseudodifferential operators on the family of mixed-norm $\alpha$-modulation spaces introduced recently  by Cleanthous and Georgiadis \cite{MR4082240}. In particular, in Section  \ref{sec:pdo} we will study pseudodifferential operators induced by symbols in the H\"ormander class
$S_{\rho}^b(\bR^n\times\bR^n)$ on mixed-norm $\alpha$-modulation spaces. For such symbols we prove in Theorem \ref{th:psibound} that
$$\sigma(x,D)\colon M^{s,\alpha}_{\vec{p},q}(\bR^n) \rightarrow
  M^{s-b,\alpha}_{\vec{p},q}(\bR^n)$$
  provided $0<\alpha\leq \rho\leq 1$,
  $\vec{p}\in (0,\infty)^n$, and $0<q<\infty$.
The case $\rho=\alpha=1$ recovers a known result that symbols in $S^b_1$ acts boundedly on the Besov spaces, see \cite{MR3573690}, but it is interesting to note that it is known that we have a strict inclusion $S^b_1\subset S^b_\rho$ for $0<\rho<1$, so a larger class of operators is covered by allowing values of $\alpha<1$. This  supports the  claim that mixed-norm $\alpha$-modulation spaces are well adapted for symbols in $S_{\rho}^b(\bR^n\times\bR^n)$.

In the non-mixed-norm case, pseudodifferential operators on $\alpha$-modulation spaces have been
considered in \cite{MR2035141, MR2035296,MR2423282,MR2292720}.  
% Pseudodifferential operators
%on $\alpha$-modulation spaces have also been studied in \cite{MR2035296,MR2423282,MR2292720}. 
Pseudodifferential operators on modulation spaces 
were first studied by Tachizawa \cite{p4},  
and later by a number of authors, see e.g.\
\cite{p2,p1,n4,MR1702232,MR1857227,p3,n5}.  In the mixed-norm setting, pseudodifferential operators on Besov and Triebel-Lizorkin spaces have been studied in \cite{MR3573690,Cleanthous2019}.

The structure of the paper is as follows. In Section
\ref{sec:mn} we introduce mixed-norm Lebesgues and $\alpha$-modulation spaces 
based on a so-called bounded admissible
partition of unity (BAPU) adapted to the mixed-norm setting.  Section
\ref{sec:mn} also introduces the maximal function estimates that will be needed to prove the main result. 
In Section \ref{sec:diff} give a precise definition of the class
$S_{\rho}^b(\bR^n\times\bR^n)$ and the associated pseudo-differential operators. We provide some boundedness results for multiplier
operators on $\alpha$-modulation spaces, and then proceed to prove the main result, Theorem \ref{th:psibound}.

\section{Mixed-norm  Spaces}\label{sec:mn}
In this section we introduce the mixed-norm Lebesgue spaces along with some needed maximal function estimates. Then we introduce mixed-norn $\alpha$-modulation spaces 
based on  so-called bounded admissible
partition of unity adapted to the mixed-norm setting. 
\subsection{Mixed-norm Lebesgue Spaces}
Let $\vec{p}=(p_1,\dots,p_n)\in(0,\infty)^n$ and $f:\mathbb{R}^n\rightarrow \mathbb{C}$. We say that $f\in L_{\vec{p}}=L_{\vec{p}}(\mathbb{R}^n)$ provided
\begin{equation}\label{Lp}
\|f\|_{\vec{p}}:=\left(\int_\mathbb{R}\cdots\left(\int_\mathbb{R}\left(\int_\mathbb{R} |f(x_1,\dots,x_n)|^{p_1} dx_1\right)^{\frac{p_2}{p_1}} dx_2\right)^{\frac{p_3}{p_2}}\cdots dx_n\right)^{\frac{1}{p_n}}<\infty,
\end{equation}
with the standard modification when $p_j=\infty$ for some $j=1,\dots,n.$ The quasi-norm $\|\cdot\|_{\vec{p}},$ is a norm when $\min(p_1,\dots,p_n)\geq 1$ and turns $(L_{\vec{p}},\|\cdot\|_{\vec{p}})$ into a Banach space. Note that when $\vec{p}=(p,\dots,p),$ then $L_{\vec{p}}$ coincides with $L_p.$ For additional properties of $L_{\vec{p}}$, see e.g.\ \cite{AI,MR370171,MR126155,Cleanthous2019}.

\subsection{Maximal operators}\label{subsec:maxop}
The maximal operator will be central to most of the estimates considered in this paper. Let $1\leq k\leq n$. We define
\begin{equation}\label{MK}
M_k f(x)=\sup\limits_{I\in I_x^k} \dfrac{1}{|I|} \int_I |f(x_1,\dots,y_k,\dots,x_n)| dy_k,\;\;f\in L^{1}_{loc}(\bR^n),
\end{equation}
where $I_x^k$ is the set of all intervals $I$ in $\mathbb{R}_{x_k}$ containing $x_k$.\

We will use extensively the following iterated maximal function:
\begin{equation}\label{Max1}
\cM_\theta f(x):=\left(M_n(\cdots(M_1|f|^\theta)\cdots)\right)^{1/\theta}(x),\;\theta>0,\;x\in\mathbb{R}^n.
\end{equation}

\begin{remark} If $Q$ is a rectangle $Q=I_1\times\dots\times I_n,$ it follows easily that for every locally integrable $f$
\begin{equation}\label{rectangle}
\int_Q |f(y)| dy\leq |Q| \cM_1 f(x)=|Q|\cM_\theta^\theta|f|^{1/\theta}(x), \ \theta>0,\ x\in\bR^n.
\end{equation}
\end{remark}

We shall need the following mixed-norm version of the maximal inequality, see \cite{MR370171,JS}:\
If $\vec{p}=(p_1,\dots,p_n)\in(0,\infty)^n$ and $0<\theta<\min(p_1,\dots,p_n)$ then there exists a constant $C$ such that
\begin{equation}\label{max}
\|\cM_\theta f\|_{\Lp}\leq C \|f\|_{\Lp}
\end{equation}

An important related estimate is a Peetre maximal function estimate, which will be one of our main tools in the sequel. 
For $ \vec{a} = (a_1,\ldots,a_n) \in\bR^n$, $\vec{b} = (b_1,\ldots,b_n) \in (0,\infty)^n$,  consider  the corresponding rectangle $R = [a_1-b_1,a_1+b_1]\times\cdots\times[a_n-b_n,a_n+b_n]$, which will be  denote $R[\vec{a},\vec{b}]$.

For every $\theta>0,\;$there exists a constant $c=c_\theta>0,$ such that for every $R>0$ and $f$ with $\supp(\hat{f})\subset c_f+R[-2,2]^n$, see \cite{MR3573690},
\begin{equation}\label{M3}
\sup_{y\in\mathbb{R}^n} \dfrac{|f(y)|}{\brac{R(x-y)}^{n/\theta}}\leq c\cM_\theta f(x),\;x\in\mathbb{R}^n.
\end{equation}
In particular, the constant $c$ is independent of the point $c_f\in \bR^n$. 
\subsection{Mixed-norm Modulation spaces}\label{sec:mod_space}
In this section we recall the definition of mixed-norm $\alpha$-modulation spaces 
as introduced by Cleanthous and Georgiadis \cite{MR4082240}. The
$\alpha$-modulation spaces form a family of smoothness spaces that contain  modulation and Besov spaces as special ``extremal'' cases. The spaces
are defined by a parameter $\alpha$, belonging  to the interval
$[0,1]$. This parameter determines a segmentation of the frequency
domain from which the spaces are built.

\begin{definition}
A countable collection $\qp$ of subsets $Q\subset \bR^n$ is called an
admissible covering of $\bR^n$ if 
\begin{itemize}
\item[i.]$\bR^n=\cup_{Q\in \qp} Q$ 
\item [ii.] There
exists $n_0<\infty$ such that 
  $\#\{Q'\in\qp:Q\cap Q'\not=\emptyset\}\leq n_0$ for all $Q\in\qp$.
\end{itemize}  
    An
  admissible covering is called an $\alpha$-covering, $0\leq
  \alpha\leq 1$, of $\bR^n$ if
 \begin{itemize}
\item[iii.] 
  $|Q|\asymp \langle x\rangle^{\alpha d}$ (uniformly) for all $x\in Q$ and for all $Q\in\qp$,
\item [iv.] There exists a constant $ K <\infty$ such that 
$$\sup_{Q\in\qp} \frac{R_Q}{r_Q}\leq K,$$
where $r_Q :=\sup\{r\in [0,\infty):\exists c_r\in\bR^n: B(c_r,r)\subseteq Q\}$
and $R_Q :=\inf\{r\in (0,\infty):\exists c_r\in\bR^n: B(c_r,r)\supseteq Q\}$, where $B(x,r)$ denotes the Euclidean ball in $\bR^n$ centered at $x$ with radius $r$.
\end{itemize}
\end{definition}

We will need a mixed-norm bounded admissible partition of unity adapted to $\alpha$-coverings. We let $\mathcal{F}(f)(\xi):=(2\pi)^{-n/2}\int_{\bR^n}
f(x)e^{- i x\cdot\xi}\,dx$, $f\in L_1(\bR^n)$,
denote the Fourier transform, and let $\hat{f}(\xi)=\mathcal{F}(f)(\xi)$.

\begin{definition}
Let  $\qp$ be an $\alpha$-covering of $\bR^n$. A corresponding
mixed-norm bounded admissible partition of unity ($\vec{p}$-BAPU) $\{\psi_Q\}_{Q\in\qp}$ is a family
of functions satisfying
\begin{itemize}
\item $\text{supp}(\psi_Q)\subset Q$
\item $\sum_{Q\in\qp} \psi_Q(\xi)=1$
\item $\sup_{Q\in \qp}|Q|^{-1}\|\chi_Q\|_{L_{\vec{\tilde{p}}}}\|\mathcal{F}^{-1}\psi_Q\|_{L_{\vec{\tilde{p}}}}<\infty$,
\end{itemize}
where $\tilde{p}_j:=\min\{1,p_1,\ldots,p_j\}$ for $j=1,2,\ldots, n$ and $\vec{\tilde{p}}:=(\tilde{p}_1,\ldots,\tilde{p}_n)$.
\end{definition}

The results in Sections
\ref{sec:diff} and \ref{sec:pdo} rely on the known fact that it is possible to construct
a smooth $\vec{p}$-BAPU with certain ``nice'' properties. We summarise the needed properties in the following proposition proved in  \cite{MR4082240}, see also \cite{Borup2006a}. Let $\brac{x}:=(1+|\,x|^2)^{1/2}$ for $x\in\bR^n$.
\begin{proposition}\label{prop:psideriv}
  For $\alpha\in [0,1)$, there exists an $\alpha$-covering of $\bR^n$
  with a corresponding $\vec{p}$-BAPU $\{\psi_k\}_{k\in\Zd}\subset \dS$ satisfying:
  \begin{itemize}
  \item[i.] $\xi_k\in Q_k$, $k\in\Zd$, where $\xi_k:= k\brac{k}^{\alpha/(1-\alpha)}$.
  \item[ii.] The following estimate holds, 
  $$|\partial^\beta \psi_k(\xi)|\leq C_\beta\langle \xi\rangle^{-|\beta|\alpha},\qquad \xi\in\bR^n,$$
  for every multi-index $\beta$ with $C_\beta$ independent of $k\in \Zd$.
  \item [iii.]  Define $\widetilde{\psi}_k(\xi) = \psi_k(|\xi_k|^\alpha\xi+\xi_k)$. Then for
every $\beta\in \bN^d$ there exists a constant
  $C_\beta$ independent of $k\in \Zd$ such that 
$$|\partial^\beta_\xi \widetilde{\psi}_k(\xi)| \leq C_\beta \chi_{B(0,r)}(\xi).$$
\item [iv.] Define
  ${\mu}_k(\xi) = \psi_k(a_k\xi)$, where
  $a_k:=\brac{\xi_k}$. Then for every $m\in \bN$ there exists a constant
  $C_m$ independent of $k$ such that 
$$|\hat{{\mu}}_k(y)| \leq C_m a_k^{(m-n)(1-\alpha)}\brac{y}^{-m},\qquad y\in\bR^n.$$
  \end{itemize}
\end{proposition}

\begin{remark}
A closer inspection of the construction presented in \cite{MR4082240} reveals that the BAPU is in fact $\vec{p}$-independent and only depends on $\alpha$ through the geometry of the $\alpha$-covering.  
\end{remark}

The case $\alpha=1$, corresponding to a dyadic-covering, is not included in Proposition \ref{prop:psideriv}, but it is known that $\vec{p}$-BAPU can easily be constructed for dyadic coverings, see  e.g.\ \cite{MR3682609}.  Using $\vec{p}$-BAPUs, it is now possible to introduce the family of mixed-norm $\alpha$-modulation spaces.

\begin{definition}\label{def:Malpha}
Let $\alpha\in [0,1], s\in \bR, \vec{p}\in (0,\infty)^n, q\in (0,\infty],$ and 
let $\qp$ be an $\alpha$-covering with associated $\vec{p}$-BAPU 
   $\{\psi_k\}_{k\in\bZ^n\backslash \{0\}}$ of the type considered in Proposition \ref{prop:psideriv}. Then we define the mixed-norm $\alpha$-{
  modulation space}, $M^{s,\alpha}_{\vec{p},q}(\bR^n)$ as the set of
  tempered distributions $f\in S^\prime(\bR^n)$ satisfying
  \begin{equation}
    \label{eq:mod}
  \|f\|_{M^{s,\alpha}_{\vec{p},q}} := \biggl( \sum_{k\in\bZ^n\backslash \{0\}}
\brac{\xi_k}^{qs} \bigl\| \cF^{-1} (\psi_k\cF f)\bigr\|_{\Lp}^q
\biggr)^{1/q}<\infty,  
  \end{equation}
 with $\{\xi_k\}_{{k\in\bZ^n\backslash \{0\}}}$ defined as in Proposition \ref{prop:psideriv}.
For $q=\infty$, we
change of the sum to $\sup_{{k\in\bZ^n\backslash \{0\}}}$.
\end{definition}

It is proved in \cite{MR4082240} that  the definition of 
$M^{s,\alpha}_{\vec{p},q}(\bR^n)$ is independent of the $\alpha$-covering and
of the BAPU, see also \cite{MR87b:46020} for the case of general decomposition space.

\section{Pseudodifferential operators on mixed-norm $\alpha$-modulation spaces}\label{sec:diff}
We now turn to the main focus of this article, the study of pseudodifferential operators on mixed-norm $\alpha$-modulation spaces. We will state and prove our main result later in this section, but let us first recall the H\"ormander class 
$S_{\rho}^b(\bR^n\times\bR^n)$, for $b\in\bR$ and $0\leq \rho\leq 1$, which 
is the family of functions $\sigma\in C^\infty(\bR^n\times\bR^n)$
satisfying
$$|\,\sigma|_{N,M}^{(b)}:=\max_{|\alpha|\leq N,|\beta|\leq
  M}\sup_{x,\xi\in\bR^n}
\brac{\xi}^{\rho|\alpha|-b}|\partial_\xi^\alpha\partial_x^\beta\sigma(\xi,x)|<\infty,$$
for all $M,N\in\bN$.  

The class $S_{\rho}^b(\bR^n\times\bR^n)$  has been studied in details in e.g.\ \cite{KG}. For $\rho<1$, we have a strict inclusion
$S_{1}^b(\bR^n\times\bR^n)\subset S_{\rho}^b(\bR^n\times\bR^n)$. An example of a symbol $\sigma\in
S_{1/2}^b(\bR\times\bR)\backslash S_{1}^b(\bR\times\bR)$ is the
symbol associated with the convolution kernel
$K(x)=e^{i/|x|}|x|^{-\gamma}$, $\gamma>0$. It can be shown that $\hat{K}(\xi)\in
S^{\gamma/2-3/4}_{1/2}(\bR^2)$,  see \cite[Chap.
VII]{MR1232192}.

We define the pseudodifferential operator $T_{\sigma}$ induced by $\sigma\in S^b_\rho$ by 
\begin{equation}
\label{dpdo}T_{\sigma}f(x):=\frac{1}{(2\pi)^{n/2}}\int_{\mathbb{R}^{n}}\sigma(x,\xi)\hat{f}(\xi)e^{ix\cdot\xi}d\xi,\ \text{for every} \ \ x\in\mathbb{R}^n, \ \ f\in\mathcal{S}(\bR^n),
\end{equation}
where $\hat{f}$ is the Fourier transform of the test function $f\in\mathcal{S}(\bR^n)$. We let $\text{Op}S^{b}_{\rho}$ denote the family of all operators induced by $S^b_\rho$. Whenever convenient, we will also use the notation $\sigma(x,D):=T_\sigma$.
 
 An important property  of  $S^b_\rho$, which we will rely on in the sequel, is the following composition result, see
e.g.\  \cite[Chap.\ 5]{KG},
 
\begin{proposition}
	\label{th:lift}
   Let $\sigma_1$ and $\sigma_2$ be
  symbols belonging to $S_{\rho}^{b_1}(\bR^n\times\bR^n)$ and
  $S_{\rho}^{b_2}(\bR^n\times\bR^n)$, respectively, for some $b_1,b_2\in \bR$. Then
  there is a symbol $\sigma\in 
  S_{\rho}^{b_1+b_2}(\bR^n\times\bR^n)$  so that $T_\sigma = T_{\sigma_1} T_{\sigma_2}$. Moreover,
  \begin{equation}
    \label{eq:Taylorsymbol}
    \sigma -\sum_{|\alpha|<N} \frac{1}{i^{|\alpha|}\alpha
      !}D_\xi^\alpha \sigma_1 \cdot D_x^\alpha \sigma_2 \in
    S_{\rho}^{b_1+b_2-N}(\bR^n\times\bR^n) ,\;\;\text{for all}\;N\in \bN.
  \end{equation}
\end{proposition}

\subsection{Fourier multipliers}
Let us first brifely a special class of pseudodifferential
operators, namely Fourier multipliers where the symbol $\sigma$ is $x$ independent. Fourier multiplies have been studied in \cite{MR4082240} and and we will just summarize the most crucial results for our study, where we will mainly rely on the Bessel potential operator. The Bessel potential  $J^b:=(I-\Delta)^{b/2}$ is defined by
$\widehat{J^bf}(\xi)=\brac{\xi}^b\hat{f}(\xi)$. It is well known that $\brac{\cdot}^b\in S^b_1$, so in particular $\brac{\cdot}^b\in S^b_\rho$ for $0<\rho\leq 1$. It also known
that for the Besov spaces we have the lifting property, $J^bB^{s}_{\vec{p},q}(\bR^n) =
B^{s-b}_{\vec{p},q}(\bR^n)$, see e.g. \cite{MR3573690}, and it was proven in \cite{MR4082240} that
$J^b$ has exactly the same lifting property when considered on
$M^{s,\alpha}_{\vec{p},q}(\bR^n)$, $0\leq \alpha\leq 1$.

\begin{proposition}\label{cor:Jb}
Let $\alpha\in [0,1]$, $s\in\bR$, $\vec{p}\in (0,\infty)^n$, $q\in(0,\infty)$. Suppose  $b\in \bR$ and let $J^b=(1-\Delta)^{b/2}$. Then  we have
  $J^bM^{s,\alpha}_{\vec{p},q}(\bR^n) =  M^{s-b,\alpha}_{\vec{p},q}(\bR^n),$
  in the sense that 
  $$\|f\|_{M^{s,\alpha}_{\vec{p},q}} \asymp
  \|J^bf\|_{M^{s-b,\alpha}_{\vec{p},q}}\qquad \text{for all}\; f\in M^{s,\alpha}_{\vec{p},q}(\bR^n).$$
\end{proposition}

\subsection{Boundedness of pseudodifferential operators}\label{sec:pdo}

We can now state and prove our main result, which we believe will provides a compelling case for the use of mixed-norm $\alpha$-modulation spaces with $\alpha<1$ as the symbol classes $S^b_\rho$ are increasing in size with  $\rho$ decreasing. 
\begin{theorem}\label{th:psibound}
  Suppose $b\in\bR$, $\alpha\in (0,1]$, $\sigma\in
  S_{\rho}^b(\bR^n\times\bR^n)$, $\alpha\leq \rho\leq 1$, $s\in \bR$,
  $\vec{p}\in (0,\infty)^n$, and
  $q\in(0,\infty)$.
  Then
  $$\sigma(x,D):M^{s,\alpha}_{\vec{p},q}(\bR^n)\rightarrow M^{s-b,\alpha}_{\vec{p},q}(\bR^n).$$ 
  Moreover, there exist $L,N>0$ (depending on
$s, \vec{p}, q$, and $\rho$) such
that the operator norm is bounded by $C|\,\sigma|_{L,N}^{(b)}$, with
$C$ a constant. 
\end{theorem}
Let us consider an example before we turn to the proof the result.
\begin{example}
Consider the
symbol associated with the convolution kernel
$K(x)=e^{i/|x|}|x|^{-\gamma}$, $\gamma>0$, $x\in \bR^2$. As mentioned earlier,  $\hat{K}(\xi)\in
S^{\gamma/2-3/4}_{1/2}(\bR^2)$. Hence, by Theorem \ref{th:psibound},
$$\hat{K}(x,D): M^{s,1/2}_{\vec{p},q}(\bR^2)\rightarrow M^{s-\gamma/2+3/4,1/2}_{\vec{p},q}(\bR^2),$$
for  $s\in \bR$,   $\vec{p}\in (0,\infty)^n$, and
  $q\in(0,\infty)$.
\end{example}

Let us now turn to the proof of Theorem \ref{th:psibound}. In the Besov space case, $\alpha=1$ [i.e.,
$M^{s,1}_{\vec{p},q}(\bR^n)=B^s_{\vec{p},q}(\bR^n)$], the proof was given by Georgiadis and the author in \cite{MR3573690}. We will therefore only consider the case $\alpha\in (0,1)$  below.

\begin{proof}[Proof of Theorem \ref{th:psibound}]

Calling on Propositions \ref{th:lift} and \ref{cor:Jb}, we have  $J^{-a}M^{s,\alpha}_{\vec{p},q}=M^{s+a,\alpha}_{\vec{p},q}$,
$\sigma(x,D) J^{a}\in \text{Op}S^{b+a}_{\rho}$, and
$J^a\sigma(x,D)\in \text{Op}S^{b+a}_{\rho}$ when $\sigma \in
S^b_{\rho}$, from which it follows that it is no restriction to assume that
$s$ is large  and $b=0$. Moreover, it suffices to prove that 
$\|\sigma(x,D)f\|_{M^{s,\alpha}_{\vec{p},q}}
\leq C \|f\|_{M^{s,\alpha}_{\vec{p},q}}$ for $f\in \dS$ since 
$\dS$ is dense in $M^{s,\alpha}_{\vec{p},q}(\bR^n)$, see \cite{MR4082240}.

Fix $f\in \dS$. First we estimate the $\Lp$-norm of $\psi_k(D)\sigma(x,D)f$.
  Notice that for any $g\in \dS$,
\begin{align}
[\psi_k(D)g](x)&=(2\pi)^{-d/2}\int_{\bR^n} e^{ix\cdot y}
  \psi_k(y)\hat{g}(y)\,dy% \nonumber\\
% &=
% (2\pi)^{-d/2}\int_{\bR^n} 
%  \check{\psi}_k(y){g}(x-y)\,dy\nonumber\\
=(2\pi)^{-d/2}\int_{\bR^n} 
  \hat{\psi}_k(y){g}(x+y)\,dy.\label{eq:est}
\end{align}
Letting $\sigma^\gamma_\eta(x,\xi):=\partial_x^\gamma\partial_\xi^\eta
\sigma(x,\xi)$ and $\sigma^\gamma:=\sigma^\gamma_0$, we obtain 
\begin{align}
\sigma(x+y,D)f(x+y)&=(2\pi)^{-n/2}\int_{\bR^n}
e^{i(x+y)\cdot \xi} \sigma(x+y,\xi)\hat{f}(\xi)\,d\xi\nonumber\\
&=
(2\pi)^{-n/2}\sum_{|\gamma|\leq K-1} \frac{y^\gamma}{\gamma!}\int_{\bR^n}
e^{i(x+y)\cdot \xi} \sigma^\gamma(x,\xi)\hat{f}(\xi)\,d\xi\nonumber\\
&+(2\pi)^{-n/2}\sum_{|\gamma|= K} K \frac{y^\gamma}{\gamma!}\int_{\bR^n}
e^{i(x+y)\cdot \xi} \int_0^1 (1-\tau)^{K-1}
\sigma^\gamma(x+\tau y,\xi)\hat{f}(\xi)\,d\tau\,d\xi\nonumber\\
&:=T(x,y)+R(x,y),\label{eq:est2}
\end{align}
where we have expanded
$\sigma(\cdot+y,\xi)$ in a Taylor series centered at $x$. We choose the order $K$ 
 such that $K\alpha>s+2(1-\alpha)(1+n)/r$, where $r:=\min\{1,q,p_1,\ldots,p_n\}$. Using \eqref{eq:est2}
in  \eqref{eq:est}, we obtain
\begin{equation}
  \label{eq:esti}
  \psi_k(D)\sigma(x,D)f(x)=(2\pi)^{-n/2}\int_{\bR^n}
\hat{\psi}_k(y)T(x,y)\,dy+
(2\pi)^{-n/2}\int_{\bR^n}
\hat{\psi}_k(y)R(x,y)\,dy.
\end{equation}
We estimate each of the two terms separately. First we consider the term
with $T(x,y)$. We have,

\begin{align}
\int_{\bR^n}
\hat{\psi}_k(y)T(x,y)\,dy&=
(2\pi)^{-n/2}\int_{\bR^n}\hat{\psi}_k(y)\sum_{|\gamma|\leq K-1} \frac{y^\gamma}{\gamma!}\
\int_{\bR^n}e^{i(x+y)\cdot \xi}
\sigma^\gamma(x,\xi)\hat{f}(\xi)\,d\xi\,dy\nonumber\\
&=(2\pi)^{-n/2}\sum_{|\gamma|\leq K-1} \frac{1}{\gamma!}
\int_{\bR^n} e^{ix\cdot\xi} \sigma^\gamma(x,\xi)\hat{f}(\xi)
\int_{\bR^n} e^{iy\cdot\xi} \hat{\psi}_k(y) y^\gamma\,dy\,d\xi\nonumber\\
&=
\sum_{|\gamma|\leq K-1} \frac{1}{\gamma!}
\int_{\bR^n}  e^{ix\cdot\xi}
\sigma^\gamma(x,\xi)\partial_\xi^\gamma
\psi_k(\xi)\hat{f}(\xi)\,d\xi.
\label{eq:est3}
\end{align}
Define $\Psi_k:= \sum_{k'} \psi_{k'}$, where
  the sum is taken over all $k'\in \bZ^n\backslash\{0\}$ with $\text{supp}(\psi_{k'})\cap \text{supp}(\psi_{k})\neq\emptyset$.  
Using the fact that $\Psi_k(\xi)=1$ on $\text{supp}(\psi_k)$, and
the relation $(\hat{f}\hat{g})^\vee = f*g$, we obtain for $\theta>0$,
\begin{align*}
\bigg| \int_{\bR^n} e^{ix\cdot\xi}
\sigma^\gamma(x,\xi)\partial_\xi^\gamma
\psi_k(\xi)\hat{f}(\xi)\,d\xi\bigg|
&
= \bigg| \int_{\bR^n} e^{ix\cdot\xi}
\bigl( \sigma^\gamma(x,\xi)\partial_\xi^\gamma
\psi_k(\xi)\bigr) \bigl( \Psi_k(\xi)\hat{f}(\xi)\bigr)\,d\xi\bigg|\\
&\leq \int_{\bR^n} 
\big|(\sigma^\gamma(x,\cdot)\partial_\xi^\gamma
\psi_k)^{\vee}(y)\big||\Psi_k(D)f(x-y)|
\,dy\\
&\leq \int_{\bR^n} \sup_{z\in\bR^n} 
\big|(\sigma^\gamma(z,\cdot)\partial_\xi^\gamma
\psi_k)^{\vee}(y)\big||\Psi_k(D)f(x-y)|
\,dy\\
&=
\int_{\bR^n} \sup_{z\in\bR^n} 
\big|(\sigma^\gamma(z,\cdot)\partial_\xi^\gamma
\psi_k)^{\vee}(y)\big|\brac{|\xi_k|^\alpha y}^{n/\theta} \\
&\qquad\qquad \times \brac{|\xi_k|^{\alpha} y}^{-n/\theta}|\Psi_k(D)f(x-y)|dy,
\end{align*}
where $\xi_k=k|k|^{\alpha/(1-\alpha)}$. Using the estimate (b) from Lemma \ref{le:tec1} below, and the Peetre maximal function estimate \eqref{M3}, we  conclude that
\begin{equation}
  \label{eq:maxineq}
  \bigg| \int_{\bR^n} e^{ix\cdot\xi}
\sigma^\gamma(x,\xi)\partial_\xi^\gamma
\psi_k(\xi)\hat{f}(\xi)\,d\xi\bigg|
\leq C |\sigma|_{L,K}^{(0)}\cM_\theta(\Psi_k(D)f)(x),
\end{equation}
with $C<\infty$ independent of $k$ and $f$, provided we choose $L>n(1+1/\theta)$. 
Hence,  we may also conclude that
\begin{equation}
  \label{eq:T}
\bigg\|\int_{\bR^n}
\hat{\psi}_k(y)T(\cdot,y)\,dy\bigg\|_{\Lp}
\leq C|\,\sigma|_{L,K}^{(0)}\|\Psi_k(D)f\|_{\Lp},
\end{equation}
provided $0<\theta<\min\{p_1,\ldots,p_n\}$. In particular, we may choose $L>n(1+1/r)$ to ensure that \eqref{eq:T}
holds.

We turn to the second term in \eqref{eq:esti}. Let
${\mu}_k(\xi)=\psi_k(a_k\xi)$, where
  $a_k:=\langle k|k|^{\alpha/(1-\alpha)}\rangle $.
First notice that
$$
\int_{\bR^n}
\hat{\psi}_k(y)R(x,y)\,dy=
\int_{\bR^n}
\hat{{\mu}}_k(y)R(x,a_k^{-1}y)\,dy.$$
We have,
\begin{align*}
\bigg|&\sum_{|\gamma|= K} \frac{a_k^{-K}}{\gamma!}\int_{\bR^n}
y^\gamma\hat{{\mu}}_k(y)
\int_{\bR^n}
e^{i(x+a_k^{-1}y)\cdot \xi} \int_0^1 (1-\tau)^{K-1}
\sigma^\gamma(x+a_k^{-1}\tau
y,\xi)\hat{f}(\xi)\,d\tau\,d\xi\,dy\bigg|\\
&\leq Ca_k^{-K}\sum_{|\gamma|= K}\int_{\bR^n}
\brac{y}^K\big|\hat{{\mu}}_k(y)\big|
\bigg|\int_0^1  (1-\tau)^{K-1}
\int_{\bR^n}
e^{i(x+a_k^{-1}y)\cdot \xi} 
\sigma^\gamma(x+a_k^{-1}\tau
y,\xi)\hat{f}(\xi)\,d\xi\,d\tau\bigg|\,dy.
\end{align*}
Using Lemma  \ref{le:tec1}  with $m=K+n+(1+n)/r$, we obtain the following estimate for the
  right-hand side for $0<\theta\leq r$, using that $0<r\leq 1$,
\begin{align*}
% with $\tilde{x}:= x-a_k^{-1}\xi_k$}
 C'a_k^{-\tilde{K}}&\sum_{|\gamma|= K}\int_{\bR^n}
\frac{\langle y\rangle^{-n-1}}{\langle y\rangle^{n/\theta}}
\sup_{z\in\bR^n}\big|[\sigma^\gamma(z,D)f](x+a_k^{-1}y)|\,dy\\
&=
 C'a_k^{-\tilde{K}}\sum_{|\gamma|= K}\int_{\bR^n}
{\langle y\rangle^{-n-1}}\sup_{z\in\bR^n}
\frac{\big|[\sigma^\gamma(z,D)f](x+a_k^{-1}y)|}{{\langle
    y\rangle^{n/\theta}}}\,dy\\
&\leq 
 C'a_k^{-\tilde{K}}\sum_{|\gamma|= K}\int_{\bR^n}
{\langle y\rangle^{-n-1}}\sup_{z,\eta\in\bR^n}
\frac{\big|[\sigma^\gamma(z,D)f](x+\eta)|}{{\langle
    a_k \eta\rangle^{n/\theta}}}\,dy\\
&\leq 
 C'a_k^{-\tilde{K}}\sum_{|\gamma|= K}\int_{\bR^n}
{\langle y\rangle^{-n-1}}\sup_{z,\eta\in\bR^n}
\frac{\big|[\sigma^\gamma(z,D)f](x+\eta)|}{{\langle
    \eta\rangle^{n/\theta}}}\,dy,\qquad\text{since $\alpha_k\geq 1$},\\
&\leq 
 C'a_k^{-\tilde{K}}\sum_{|\gamma|= K}
\sup_{z,\eta\in\bR^n}
\frac{\big|[\sigma^\gamma(z,D)f](x+\eta)|}{{\langle
    \eta\rangle^{n/\theta}}},
\end{align*}
where $\tilde{K}=K\alpha-(1+n)(1-\alpha)/r \geq s+\frac{n+1}{q}(1-\alpha)$, since $q\geq r$.
Now,
\begin{align*}
\bigg(\sum_{k\in\Zd} a_k^{sq}\bigg\|&\int_{\bR^n}
\hat{\psi}_k(y)R(x,y)\,dy\bigg\|_{\Lpx}^q\bigg)^{1/q}\\
&\leq 
\bigg\{C\sum_{k\in\Zd}
a_k^{(s- \tilde{K})q}\bigg(\sum_{|\gamma|= K}\biggl\|
\sup_{\eta,z\in\bR^n}
\frac{\big|[\sigma^\gamma(z,D)f](x+\eta)|}{{\langle
    \eta\rangle^{n/\theta}}}\bigg\|_{\Lpx} \bigg)^q\bigg\}^{1/q}.\\
\intertext{We notice that $L^q:=C\sum_{k\in\Zd}
a_k^{(s- \tilde{K})q}\leq C\sum_{k\in\Zd}
|k|^{-n-1}$ is finite. Based on this observation,  
recalling that $r=\min\{1,q,p_1,\ldots,p_n\}$, we use the equivalence of $\ell^\tau$-norms on finite dimensional spaces to estimate the right-hand side by,}
 L
\bigg(\sum_{|\gamma|= K}&\bigg\|
\sup_{\eta,z\in\bR^n}
\frac{\big|[\sigma^\gamma(z,D)f](x+\eta)|}{{\langle
    \eta\rangle^{n/\theta}}}\bigg\|_{\Lpx}^r\bigg)^{1/r}\\
&= L
\bigg(\sum_{|\gamma|= K}\bigg\|
\sup_{\eta,z\in\bR^n}
\frac{\big|[\sigma^\gamma(z,D)\sum_{k\in\Zd} \psi_k(D)f](x+\eta)|}{{\langle
     \eta\rangle^{n/\theta }}}\bigg\|_{\Lpx}^r\bigg)^{1/r}\\
&\leq  L
\bigg(\sum_{|\gamma|= K}\sum_{k\in\Zd}\bigg\|
\sup_{\eta,z\in\bR^n}
\frac{\big|[\sigma^\gamma(z,D) \psi_k(D)f](x+\eta)|}{{\langle
    \eta\rangle^{n/\theta}}}\bigg\|_{\Lpx}^r \bigg)^{1/r}.
\end{align*}
We now focus on the individual term $A_k:=\big|[\sigma^\gamma(z,D)
\psi_k(D)f](x+\eta)|$. Put $f_k(x):=[\Psi_k(D)f](x)$, with $\Psi_k$ defined as above. We have
\begin{align*}
A_k&=
\bigg|\int_{\bR^n}
(\sigma^\gamma(z,\xi)\psi_k(\xi))^\vee(x+\eta-y)f_k(y)\,dy\bigg|\\
&\leq\int_{\bR^n}|
(\sigma^\gamma(z,\xi)\psi_k(\xi))^\vee(x+\eta-y)||f_k(y)|\,dy\\
&\leq
\sup_{u\in \bR^n} \frac{|f_k(u)|}{\langle x-u\rangle^{n/\theta}}
\int_{\bR^n} |
(\sigma^\gamma(z,\xi)\psi_k(\xi))^\vee(x+\eta-y)|\langle x-y\rangle^{n/\theta}\,dy.
\end{align*}
Now, $\langle x-y\rangle^{n/\theta}\leq
c\langle x-y+\eta\rangle^{n/\theta} \langle \eta\rangle^{n/\theta}$, so
\begin{align*}
\sup_{z,\eta\in \bR^n}\frac{A_k}{\langle \eta\rangle^{n/\theta}}
&\leq C
\sup_{\eta\in \bR^n} \frac{|f_k(x-\eta)|}{\langle
  \eta\rangle^{n/\theta}}
\sup_{z\in \bR^n}
\int_{\bR^n}
|(\sigma^\gamma(z,\xi)\psi_k(\xi))^\vee(u)|\langle u
\rangle^{n/\theta }\,du\\
&\leq C'
\sup_{\eta\in \bR^n} \frac{|f_k(x-\eta)|}{\langle
  \eta\rangle^{n/\theta}} |\,\sigma|^{(0)}_{L,K},
\end{align*}
provided $L>n+n/\theta$, where we have used Lemma \ref{le:tec1}.
Hence,
\begin{align*}
  \sum_{|\gamma|= K}\sum_{k\in\Zd}\bigg\|
  \sup_{\eta,z\in\bR^n}&
  \frac{\big|[\sigma^\gamma(z,D) \psi_k(D)f](x+\eta)|}{{\langle
      \eta\rangle^{n/\theta }}}\bigg\|_{\Lpx}^r\\
  &=\sum_{|\gamma|= K}\sum_{k\in\Zd}
  a_k^{nr/\theta}
  \bigg\|
  \sup_{\eta,z\in\bR^n}
  \frac{\big|[\sigma^\gamma(z,D)
    \psi_k(D)f](x+\eta)|}{a_k^{ n/\theta}
      \langle
      \eta\rangle^{n/\theta}}\bigg\|_{\Lpx}^r\\
  &\leq C'
  (|\,\sigma|^{(0)}_{L,K})^r \sum_{k\in\Zd} a_k^{nr/\theta }
  \biggl\| \sup_{\eta\in \bR^n}
  \frac{|f_k(x-\eta)|}{\langle a_k
    \eta\rangle^{n/\theta}}\biggr\|_{\Lpx}^r. 
\end{align*}
 We now use the Peetre maximal estimate,
$$\sup_{z\in\bR^n}\frac{|f_k(x-z)|}{\langle a_k z\rangle^{n/\theta }}\leq C \cM_\theta(f_k)(x),$$
and we may apply the $L_{\vec{p}}$-norms, using the
maximal inequality, to obtain
\begin{align*}
\bigg\|\sup_{z\in\bR^n}\frac{|f_k(x-z)|}{\langle
  a_k z\rangle^{n/\theta}}\bigg\|_{\Lp}^r&\leq
C'\|f_k\|_{\Lp}^r,
\end{align*}
provided $0<\theta<\min\{p_1,\ldots,p_n\}$.
Putting these estimates together yields,
\begin{align*}
%  \label{eq:rema}
  \sum_{|\gamma|= K}\sum_{k\in\Zd}\bigg\|
\sup_{\eta,z\in\bR^n}&
\frac{\big|[\sigma^\gamma(z,D) \psi_k(D)f](x+\eta)|}{{\langle
   \eta\rangle^{n/\theta}}}\bigg\|_{\Lpx}^r\\&
   \leq C''(|\,\sigma|^{(0)}_{L,K})^r
\sum_{k\in\Zd} a_k^{nr/\theta}\|\Psi_k(D)f\|_{\Lp}^r,
\end{align*}
provided $L>n+n/\theta$,
and consequently
\begin{align*}
%  \label{eq:rema2}
\bigg(\sum_{k\in\Zd} a_k^{sq}\bigg\|\int_{\bR^n}
\hat{\psi}_k(y)R(x,y)\,dy\bigg\|_{L_p}^q\bigg)^{1/q}
&\leq C''|\,\sigma|^{(0)}_{L,K}\bigg(
\sum_{k\in\Zd} a_k^{r n\theta }\|\Psi_k(D)f\|_{\Lp}^r\bigg)^{1/r}.\\
&= C''|\,\sigma|^{(0)}_{L,K}\bigg(
\sum_{k\in\Zd} a_k^{r n\theta-sr }\cdot a_k^{sr} \|\Psi_k(D)f\|_{\Lp}^r\bigg)^{1/r}.\\
&\leq C'''|\,\sigma|^{(0)}_{L,K}
\biggl\| a_k^s\|\Psi_k(D)f\|_{\Lp}\biggr\|_{\ell_q},
\end{align*}
where for the last estimate, we used H\"older's inequality with parameters $q/r$ and $q/(q-r)$ and the fact that  $s>n(1+\theta)/r$, where we also notice that $n(1+\theta)/r<3n/r$ since $0<\theta<2$.
Finally, we can put all the estimates together to close the case $b=0$ and
$s>3n/r$. We have, with $L>n+n/r$,
\begin{align*}
\|&\sigma(x,D)f\|_{M^{s,\alpha}_{\vec{p},q}}\\
&\asymp\big\|a_k^s\|\psi_k(D)\sigma(x,D)f\|_{\Lp}\big\|_{\ell^q}\\
&\leq
C\bigg\{
\bigg\|a_k^s\bigg\|\int_{\bR^n}
\hat{\psi}_k(y)T(x,y)\,dy\bigg\|_{\Lpx}\bigg\|_{\ell_q}
+\bigg\|a_k^s\bigg\|\int_{\bR^n}
\hat{\psi}_k(y)R(x,y)\,dy\bigg\|_{\Lpx}\bigg\|_{\ell_q}\bigg\}\\
&\leq C'\bigg(
|\,\sigma|^{(0)}_{L,K}\bigg\| a_k^s\big\|\Psi_k(D)f\big\|_{\Lp}\bigg\|_{\ell_q}
+|\,\sigma|^{(0)}_{L,K}\bigg\|a_k^s\big\|\Psi_k(D)f\big\|_{\Lp}\bigg\|_{\ell_q}\bigg)\\
&\leq C'' |\,\sigma|^{(0)}_{L,K}\|f\|_{M^{s,\alpha}_{\vec{p},q}}.
\end{align*}
This concludes the proof of the theorem.
\end{proof}

%\begin{remark}
%A closer examination of the arguments used in the proof reveals  that there exist $M,N>0$ (depending on
%$s, q$, and $\rho$) such
%that
%the norm of the operator $$\sigma(x,D)\colon M^{s,\alpha}_{p,q}(\bR^n)\rightarrow
%M^{s-b,\alpha}_{p,q}(\bR^n)$$ is bounded by $C|\,\sigma|_{M,N}^{(b)}$, with
%$C$ a constant. 
%\end{remark}
The following technical lemma has been used in the proof of Theorem \ref{th:psibound}.
\begin{lemma}\label{le:tec1}
Let $\alpha\in [0,1)$ and let $\{\psi_k\}_{k\in \bZ^n\backslash\{0\}}$ be the $\vec{p}$-BAPU from Proposition \ref{prop:psideriv}, depending only on $\alpha$. Suppose $\sigma \in S^0_{\rho,0}$, $\alpha\leq \rho\leq 1$. Then for
  $|\gamma|\leq K$ and $m\geq 0$, we have 
 \begin{enumerate}
\item[(a)]
For
 $|\gamma|,|\nu|\leq K$ and $J\in\bN$ 
there exists a constant $C:=C(K,J)$ such that
 $$M(x):= \sup_{z\in\bR^n}\big|
 (\partial_x^\gamma\sigma(z,\cdot)\partial_\xi^\nu
 \psi_k)^{\vee}(x)\big|\leq
 C
 |\,\sigma|_{J,K}^{(0)} |k|^{\alpha n/(1-\alpha)}\brac{|k|^{\alpha/(1-\alpha)} x}^{-J},$$
 for $x\in\bR^n, k\in\bN$.

 \item[(b)] For
  $|\gamma|,|\nu|\leq K$ and $m\geq 0$ there exists a constant 
$C':=C'(K,m)$,
  such that
  $$I:=\int_{\bR^n}\sup_{z\in\bR^n} \big|
  (\partial_x^\gamma\sigma(z,\cdot)\partial_\xi^\nu
  \psi_k)^{\vee}(x)\big|\langle |k|^{\alpha/(1-\alpha)}x \rangle^m\,dx\leq
  C' |\,\sigma|_{M,K}^{(0)},\qquad k\in\bN,$$
  for any $M\in \bN$ satisfying $M>m+n$.
 \end{enumerate}
\end{lemma}
\begin{proof}
  First we prove (a). Let
  $\sigma^\gamma_\eta(x,\xi):=\partial_x^\gamma\partial_\xi^\eta
  \sigma(x,\xi)$ and $\sigma^\gamma:=\sigma^\gamma_0$.  We have the equality
  $$M(x)=  \sup_{z\in\bR^n} 
  \bigg|\int_{\bR^n} e^{ix\cdot \xi}\sigma^\gamma(z,\xi)\partial_\xi^\nu
  \psi_k(\xi)\,d\xi\bigg|.$$
  % Let $\widetilde{\psi}_k(\xi):=
  % \psi_k(T_k\xi)$, and
 Let $T_k=|\xi_k|^\alpha+\xi_k$, where $\xi_k=k|k|^{\alpha/(1-\alpha)}$. Then a substitution yields    
  \begin{align}\label{eq:m}
   M(x)= |\xi_k|^{n\alpha}\sup_{z\in\bR^n} 
    \bigg|\int_{\bR^n} e^{i |\xi_k|^{\alpha} x\cdot \xi}\sigma^\gamma(z,T_k\xi)\partial_\xi^\nu
    \psi_k(T_k\xi)\,d\xi\bigg|.
  \end{align}
Fix $J>1$. We use the well-known estimate 
$\langle x\rangle^J|\hat{g}(x)|\leq C_J \sum_{|\beta|\leq
  J}\|\partial^\beta{g}\|_{L_1}$, for 
some finite constant $C_J$. We apply the estimate to \eqref{eq:m} to obtain
  $$M(|\xi_k|^{-\alpha} x)\leq C_J
|\xi_k|^{n\alpha}
  \sup_{z\in\bR^n} 
  \sum_{|\beta|\leq J}
  \int_{\bR^n} \Bigl| \partial_\xi^\beta \Bigl[
  \sigma^\gamma (z,T_k\xi)
  \partial_\xi^{\nu}
  \psi_k (T_k\xi)\Bigr]\Bigr|\,d\xi
  \langle x\rangle^{-J},$$
  which by Leibniz's rule provides the bound
  \begin{align}
    \label{eq:Iest}
M(|\xi_k|^{-\alpha} x)\leq
C'|\xi_k|^{n\alpha}
      \mathop{\sum_{|\beta|\leq J}}_{0\leq \eta\leq \beta} 
    \sup_{z\in\bR^n}
    \int_{\bR^n} |\xi_k|^{\alpha|\eta|}|\partial_\xi^\eta(\sigma^\gamma(z,T_k\xi))|
    |\partial_\xi^{\beta-\eta}(\partial_\xi^{\nu}
    \psi_k(T_k\xi))|\,d\xi  \langle
    x\rangle^{-J}.
  \end{align}
From Proposition \ref{prop:psideriv}, we have 
 \begin{equation}
    \label{eq:partest}
  |\partial_\xi^{\beta-\eta}(\partial_\xi^{\nu} \psi_k(T_k\xi))|
  \leq C\chi_{Q}(\xi),  
  \end{equation}  
with $C:=C(\nu,\beta,\eta)$.   
 %
%  Let us take a closer look at $|\partial_\xi^{\beta-\eta}(\partial_\xi^{\alpha}
%    \psi_k(T_k\xi))|$. Put $\mu_k(\xi) =\psi_k(T_k\xi)$. It is easy to
%    verify that $|\partial_\xi^\beta \mu_k(\xi)|\leq C_\beta
%    \chi_{Q}(\xi)$ where $Q\subset \bR^n$ is a fixed compact set
%    independent of $k$, see Equation \eqref{eq:babound}. Notice that the
%    chain rule 
%    yields 
%$$\partial_\xi^{\alpha} \psi_k(T_k\xi) =\sum_{|\beta|=|\alpha|} b_\beta
%\partial_\xi^\beta \mu_k(\xi),$$
%where $b_\beta$ are monomials of degree $|\beta|$ in the entries of
%$\delta_{h(\xi_k)^{-1}}$. Now, since $h(\xi)\geq \epsilon_0>0$, each entry is
%uniformly bounded, and we obtain
 % \begin{equation}
%    \label{eq:partest}
%  |\partial_\xi^{\beta-\eta}(\partial_\xi^{\alpha} \psi_k(T_k\xi))|
%  \leq C\chi_{Q}(\xi).  
%  \end{equation}  
We also notice that for
$\xi\in Q$,
\begin{equation}\label{eq:partest2}
	|\partial_\xi^\eta(\sigma^\gamma(z,T_k\xi))|\leq
|\,\sigma|_{|\eta|,K}^{(0)} \brac{|\xi_k|^\alpha\xi+\xi_k}^{-\rho|\eta|}
\leq C|\,\sigma|_{|\eta|,K}^{(0)} \brac{\xi_k}^{-\rho|\eta|}.
\end{equation}
 
Now,  by assumption $\alpha\leq \rho$, so using the estimates \eqref{eq:partest} and \eqref{eq:partest2} in
  \eqref{eq:Iest}, we obtain
  \begin{equation*}
M(|\xi_k|^{-\alpha} x)\leq
C''|\xi_k|^{n\alpha} 
     \mathop{\sum_{|\beta|\leq J}}_{0\leq \eta\leq \beta}
    |\,\sigma|_{J,K}^{(0)}
    \int_{\bR^n} \chi_{Q}(\xi)\,d\xi  \langle
    x\rangle^{-J}
    \leq C''' |\xi_k|^{\alpha n}\cdot |\,\sigma|_{J,K}^{(0)}\langle
    x\rangle^{-J},
  \end{equation*}
  which proves (a), since $|\xi_k|=|k|^{1/(1-\alpha)}$. 

Let us turn to (b). Pick $J>m+n$ in (a). We have
\begin{align*}
  I&=\int_{\bR^n} M(x)\brac{|k|^{\alpha/(1-\alpha)}x}^m\, dx\\
  &\leq C' 
  |\,\sigma|_{J,K}^{(0)} |k|^{\alpha n/(1-\alpha)}
 \int_{\bR^n}
  \brac{|k|^{\alpha/(1-\alpha)} x}^{-J}\brac{|k|^{\alpha/(1-\alpha)}x}^m \,dx\\
  &= C' |\,\sigma|_{J,K}^{(0)} \int_{\bR^n}
  \brac{x}^{-J}\brac{ x}^{m}\, dx\\
  &\leq 
\tilde{C}|\,\sigma|_{J,K}^{(0)},
\end{align*}
where we made a change of variable in the integral and used $J>n+m$, which of course implies that $m-J<-n$. This concludes the proof. 
\end{proof}

\section{Hypoelliptic pseudodifferential operators}\label{sec:elliptic}
In this final section we consider an application of the result in the
previous section to hypoelliptic pseudodifferential operators based on standard machinery,  see e.g.\ 
\cite{MR781536}. 
Let us introduce some notation. Let
$$S^\infty_{\rho} := \bigcup_{b\in \bR}
S^{b}_{\rho},\quad \text{and}\quad S^{-\infty}_{\rho} := \bigcap_{b\in \bR}
S^{b}_{\rho}.$$
Assume that $b_0,b\in\bR$ such that $b_0\leq b$. An element
$\sigma\in S^b_{\rho}(\bR^n\times\bR^n)$ is called {\em hypoelliptic} with
parameters $b_0$ and $b$ if there are positive constants $c$ and $a$ such that
$$a\brac{\xi}^{b_0}\leq |\sigma(x,\xi)|,\qquad \brac{\xi}\geq c,$$
and
$$|\partial_\xi^\alpha\partial_x^\beta\sigma(x,\xi)|\leq
C_{\alpha,\beta}
|\sigma(x,\xi)|\brac{\xi}^{-\rho|\alpha|},\qquad
\brac{\xi}\geq c.$$
Let $HS^{b,b_0}_{\rho}(\bR^n\times\bR^n)$ the family of all
such symbols.
The following result is  well-know, see \cite[Theorem
22.1.3]{MR781536}.
\begin{theorem}\label{th:folland}
  Suppose $\sigma \in HS^{b,b_0}_{\rho}$, with $0<\rho\leq 1$.
  Then there exists $\tau \in HS^{-b_0,-b}_{\rho}$ such that
  $I-\sigma(x,D)\tau(x,D)$ and $I-\tau(x,D)\sigma(x,D)$ are both in
  $\text{Op}(S^{-\infty}_{\rho})$.
\end{theorem}

Let $M^{-\infty,\alpha}_{\vec{p},q}(\bR^n) = \cup_{s\in \bR} M^{s,\alpha}_{\vec{p},q}(\bR^n)$.
Using Theorem \ref{th:folland} and the result from the previous
section we have

\begin{theorem}
  Suppose $\sigma \in HS^{b,b_0}_{\rho}$, with 
  $\rho\geq \alpha>0$, and $f\in M^{-\infty,\alpha}_{p,q}(\bR^n)$. If
  $\sigma(\cdot,D)f \in M^{s,\alpha}_{\vec{p},q}(\bR^n)$ for some $s\in
  \bR$, then $f\in M^{s+b_0,\alpha}_{\vec{p},q}(\bR^n)$.
\end{theorem}

\begin{proof}
  Let $S=\sigma(\cdot,D)$, and let $T=\tau(\cdot,D)$ be as in Theorem
  \ref{th:folland}. Notice that $f=T(Sf)+(I-TS)f$. By Theorem \ref{th:psibound}, $T$
  maps $M^{s,\alpha}_{\vec{p},q}(\bR^n)$ to $M^{s+b,\alpha}_{\vec{p},q}(\bR^n)$ and $(I-TS)$
  maps $M^{-\infty,\alpha}_{\vec{p},q}(\bR^n)$ to $M^{s+b,\alpha}_{\vec{p},q}(\bR^n)$.
\end{proof}

The following example will conclude the paper.

\begin{example}\label{ex:heat}
Consider the heat operator $L$ given by
$$L(u):=\frac{\partial u}{\partial t}-\sum_{j=1}^d \frac{\partial^2
  u}{\partial x_j^2}.$$
The symbol of $L$ is given by
$$l(\tau,\xi)=(i \tau+|\xi|^2),\qquad (\tau,\xi)\in \bR\times\bR^n,$$
and one can easily verify that $l\in HS^{2,1}_{1}$.
We consider an approximate inverse $P$ to $L$ with symbol
 $$a(\tau,\xi)=(i \tau+|\xi|^2)^{-1}\eta(\tau,\xi),\qquad (\tau,\xi)\in \bR\times\bR^n,$$
where $\eta$ is a smooth cut-off function that vanishes near the
origin and equals 1 for large $(\tau,\xi)$. It is easy to check  that
$a(\tau,\xi)\in
HS^{-1,-2}_{1}(\bR^{n+1}\times\bR^{n+1})$.
Hence, if $u\in
M^{-\infty,\alpha}_{\vec{p},q}(\bR^{n+1})$, $\vec{p}\in (0,\infty)^n$, $0<q<\infty$, $\alpha\in (0,1]$, and $P(u)\in
M^{s,\alpha}_{\vec{p},q}(\bR^{n+1})$, then $u\in M^{s-1,\alpha}_{\vec{p},q}(\bR^{n+1})$.
\end{example}

%\bibliographystyle{abbrv}

%\bibliography{new_pdo}

\end{document}